\documentclass[12pt,reqno]{amsart}

\newcommand\version{July 11, 2022}

\usepackage{amsmath,amsfonts,amsthm,amssymb,amsxtra}
\usepackage{bbm} 

\theoremstyle{definition}

\theoremstyle{remark}

\renewcommand{\epsilon}{\varepsilon}

\newcommand{\N}{\mathbb{N}}

\renewcommand{\phi}{\varphi}
\newcommand{\R}{\mathbb{R}}

\DeclareMathOperator{\spec}{spec}

\DeclareMathOperator{\Tr}{Tr}

\begin{document}

\title[The work of Elliott Lieb --- \version]{The work of Elliott Lieb}

\author{Rupert L. Frank}
\address[Rupert L. Frank]{Mathe\-matisches Institut, Ludwig-Maximilans Universit\"at M\"unchen, The\-resienstr.~39, 80333 M\"unchen, Germany, and Munich Center for Quantum Science and Technology, Schel\-ling\-str.~4, 80799 M\"unchen, Germany, and Mathematics 253-37, Caltech, Pasa\-de\-na, CA 91125, USA}
\email{r.frank@lmu.de}

\begin{abstract}
	On the occasion of Elliott Lieb being awarded the Gauss Prize 2022, we give a non-technical overview over some of his seminal works in mathematical physics. We emphasize, in particular, his work on Coulomb many-body systems and functional inequalities.
\end{abstract}

\renewcommand{\thefootnote}{${}$} \footnotetext{\copyright\, 2022 by the author. This paper may be reproduced, in its entirety, for non-commercial purposes.\\
	Partial support through US National Science Foundation grants DMS-1954995, as well as through the Deutsche Forschungsgemeinschaft (DFG, German Research Foundation) through Germany’s Excellence Strategy EXC-2111-390814868 is acknowledged.\\
The author is grateful to M. Lewin for helpful remarks.}

\maketitle

Elliott Lieb is awarded the Gauss Prize 2022 ``for deep mathematical contributions of exceptional breadth which have shaped the fields of quantum mechanics, statistical mechanics, computational chemistry, and quantum information theory". 

It is my great pleasure to congratulate Elliott on this honor. In the following pages I will try to give a non-technical overview over some of his seminal works.

Lieb is a mathematical physicist. This is a field that lies at the boundary between physics and mathematics and that Freeman Dyson [From Eros to Gaia, pp.\ 164--165] has described as follows:
\begin{quote}
	``Mathematical physics is the discipline of people who try to reach a deep understanding of physical phenomena by following the rigorous style and method of mathematics."
\end{quote}
As mentioned in the citation, Lieb has made ground-breaking contributions to both mathematics and physics. In this connection it is worthwhile to mention that about half a year ago, Lieb was awarded the 2022 APS Medal for Exceptional Achievement in Research, the highest honor of the American Physical Society.

A distinctive feature of Lieb's work is its timelessness. Independently of fashions and trends, he has worked, and continues to work, on problems that he considers deep and fundamental. Sometimes it is decades later that the full potential of his ideas is understood. Prime examples are the Lieb--Robinson bounds, proved in 1972, and the strong subadditivity of entropy, proved in 1973, that have both come to play a key role in quantum information theory in the 21st century.

Lieb's publication list contains, at the time of this writing, 404 items with the first one dating back to 1955 and the most recent one just about to appear. Four volumes of his selected works have been published so far \cite{selecta1,selecta2,selecta3,selecta4} and, on the occasion of his 90th birthday, a 1300+ page collection of articles was edited \cite{birthday}, where the contributors explain the content and the ramifications of Lieb's work.

It is impossible to give an overview over this monumental body of work. I have chosen here two main areas of Lieb's research, namely quantum many-body Coulomb systems and functional inequalities, and omitted all others, except for a brief mention of some in the last section. This selection is influenced by my predilections and my ignorance, and I ask the readers' indulgence for all the omissions.


\section{Quantum Coulomb systems}

One of the recurring themes in Lieb's research since the early 1970s is the analysis of continuous quantum many-body systems and, in particular, of systems of particles interacting through Coulomb forces. Here one neglects other forces, but this description is appropriate for ordinary matter and much of the world relevant to everyday life. As a very readable account of Lieb's research program in this area we recommend his Gibbs lecture in 1989 \cite{gibbs}, which contains a much more detailed description than what we provide here.

We consider a system consisting of $N$ quantum electrons and $K$ classical nuclei in $\R^3$. The latter are fixed at positions $R_1,\ldots,R_K\in\R^3$ and have charges $Z_1,\ldots,Z_K\in (0,\infty)$ (in units where the electron charge is $-1$). The properties of such a system are described by the Hamiltonian
\begin{equation}
	\label{eq:hamiltonian}
	H :=
	\sum_{n=1}^N (-\Delta_n) - \sum_{n=1}^N\sum_{k=1}^K \frac{Z_k}{|x_n-R_k|} + \sum_{1\leq n< m\leq N} \frac{1}{|x_n-x_m|} + \sum_{1\leq k<\ell\leq K} \frac{Z_k Z_\ell}{|R_k-R_\ell|} \,,
\end{equation}
acting on functions in $\R^{3N}$. Here we use coordinates $x=(x_1,\ldots,x_N)\in\R^{3N}$. The four sums in the definition of $H$ correspond, respectively, to the kinetic energy of the electrons, the electron-nucleus attraction, the electron-electron repulsion and the nucleus-nucleus repulsion.

This Hamiltonian can be realized as an unbounded, selfadjoint operator in the Hilbert space consisting of all antisymmetric functions in $L^2(\R^{3N})$, that is, square-integrable functions $\psi$ satisfying $\psi(\ldots,x_n,\ldots,x_m,\ldots)=- \psi(\ldots,x_m,\ldots,x_n,\ldots)$ for all $n\neq m$. The antisymmetry reflects the Pauli exclusion principle. One should also take the electron spin into account, but mathematically this only leads to minor changes and will be ignored in what follows.

The ground state energy is, by definition,
$$
\inf\spec H \,.
$$
It is the infimum of $\langle\psi,H\psi\rangle$ with respect to all normalized, antisymmetric $\psi$.

The fundamental feature of the problem of analyzing the ground state energy of $H$ (and of quantum many-body systems in general) is the huge number of dimensions of the underlying Hilbert space. This makes a numerical computation virtually impossible and for quantitative results one typically has to rely on approximate theories that are numerically more tractable. This situation underlines the importance of analytical studies about the full Schr\"odinger problem and also about its relation to approximate theories. Lieb has made fundamental contributions to this problem, as we will review in the remainder of this section.


\subsection*{Stability of matter}

The problem of \emph{stability of matter} consists in showing that for every $z>0$ there is a constant $C$ such that, for all $N,K\in\N$, $R_1,\ldots, R_K\in\R^3$ and $Z_1,\ldots, Z_K\in[0,z]$, one has
\begin{equation}
	\label{eq:som}
	\inf\spec H \geq - C\, (N+K) \,.
\end{equation}
This inequality says that, despite the fact that the number of interactions grows \emph{quadratically} in the total number of particles, the ground state energy only behaves \emph{linearly}. This is fundamental for the existence of matter as we know it.

We emphasize that the stability of matter depends on the fact that electrons are fermions, that is, on the fact that $H$ is considered only on the subspace of antisymmetric functions in $L^2(\R^{3N})$ and not on the full space $L^2(\R^{3N})$. It is, in fact, wrong on the latter, bigger space, as we shall see.

The first proof of stability of matter was achieved by Dyson and Lenard in 1967. Lieb and Thirring gave a new proof in 1975 \cite{lt}. The latter proof gives a much more realistic value for the constant $C$ in \eqref{eq:som}, namely about $C\approx 5$ instead of Dyson and Lenard's $C\approx 10^{14}$. It also helps to clarify the reason for the validity of stability of matter on a conceptual level. Here is what Dyson writes in the preface of the Selecta of Elliott Lieb \cite{selecta1}:
\begin{quote}
``Our proof was so complicated and so unilluminating that it stimulated Lieb and Thirring to find the first decent proof [...]. Why was our proof so bad and why was theirs so good? The reason is simple. Lenard and I began with mathematical tricks and hacked our way through a forest of inequalities without any physical understanding. Lieb and Thirring began with physical understanding and went on to find the appropriate mathematical language to make their understanding rigorous. Our proof was a dead end. Theirs was a gateway to the new world of ideas [...].''
\end{quote}

The fundamental mechanism of the Lieb--Thirring proof of stability of matter is the fact that atoms do not bind in an approximate model of a Coulomb system known as Thomas--Fermi theory. (Here ``no binding" means, mathematically, that the ground state energy of a molecule is simple the sum of the energies of the individual atoms.) Therefore the goal of Lieb and Thirring was to show that this approximate theory provides, after an appropriate adjustment of constants, a rigorous lower bound to the Schr\"odinger ground state energy. Later in this section we will discuss both this Thomas--Fermi theory and a functional inequality (known as Lieb--Thirring inequality) that leads to the claimed lower bound.

The above mentioned work by Lieb and Thirring from 1975 was the starting point of Lieb's thorough investigation of the problem of stability of matter, a topic to which he would return repeatedly for several decades and with many coworkers. Notable are, among other things, a proof of stability of matter in the presence of magnetic fields that covers the physical value of the fine structure constant \cite{sommag}. An introduction to this area is provided in his book with Seiringer ``The Stability of Matter in Quantum Mechanics” \cite{sombook}.


\subsection*{Existence of the thermodynamic limit for real matter with Coulomb forces}

Lieb's first result on Coulomb many-body systems, even before his work on stability of matter, settled an open problem in the foundations of statistical mechanics. Namely, in 1969 he and Lebowitz proved the existence of the thermodynamic limit for real matter with Coulomb forces \cite{thermo}.

Here one considers a large number of particles (electrons and nuclei, for instance) confined to an open set $\Omega$ in $\R^3$. One considers the limit where $\Omega$ tends to $\R^3$ (in a sense to be made precise) and where the densities of the different particles (that is, the number of particles divided by the volume of $\Omega$) tend to given constants. The particles interact through Coulomb forces, as in \eqref{eq:hamiltonian}. Moreover, one typically discusses this question at a given positive temperature. The theorem of Lieb and Lebowitz states that the limit of the corresponding free energy exists, is independent of the shape of the approximating domains $\Omega$ and has appropriate convexity and concavity properties.

The Lieb--Lebowitz theorem uses the stability of matter theorem as an ingredient. In the proof of the existence of the thermodynamic limit the main concern is the slow decay of $|x|^{-1}$ as $|x|\to\infty$. Indeed, the existence of the thermodynamic limit had been known in the case of short-range interactions. In the long range case charge neutrality is an essential input. Lieb and Lebowitz exploit the electrostatic screening very originally via Newton’s theorem. In this way they are led to the geometric problem of how to efficiently pack large balls by smaller balls, which they solve by their ``Cheese Theorem''.


\subsection*{Thomas--Fermi theory and Density Functional Theory}

In our discussion of the Lieb--Thirring proof of stability of matter we have already mentioned the Thomas--Fermi model of a Coulomb system. This was proposed independently by Thomas and Fermi in 1927, remarkably soon after Schr\"odinger had introduced his theory. In the variational approach to Thomas--Fermi theory one starts from the energy functional, defined for nonnegative function $\rho$ on $\R^3$ by
$$
\mathcal E[\rho] = \gamma^{\rm TF} \int_{\R^3} \rho(x)^\frac{5}{3}\,dx - \sum_{k=1}^K \int_{\R^3} \frac{Z_k\ \rho(x)}{|x-R_k|}\,dx + D[\rho] + \sum_{1\leq k<\ell\leq K} \frac{Z_k Z_\ell}{|R_k-R_\ell|}
$$
with
$$
D[\rho] := \frac12 \iint_{\R^3\times\R^3} \frac{\rho(x)\,\rho(y)}{|x-y|}\,dx\,dy \,.
$$
The function $\rho$ describes the distribution of electrons and its integral has the meaning of the total number of electrons. The four terms in the definition of $\mathcal E$ correspond, respectively, to the kinetic energy of the electrons, the electron-nucleus attraction, the electron-electron repulsion and the nucleus-nucleus repulsion. In the next subsection we will briefly describe how one can arrive at the $\rho^\frac53$ approximation to the kinetic energy. This approximation also leads to a certain value for the constant  $\gamma^{\rm TF}>0$. 

The ground state energy in Thomas--Fermi theory is
$$
\inf\left\{ \mathcal E[\rho] :\ \rho\geq 0 \,,\ \int_{\R^3} \rho\,dx = N \right\} \,.
$$

We emphasize that, in contrast to Schr\"odinger theory, Thomas--Fermi theory is a nonlinear theory. The important simplifying feature is that in Thomas--Fermi theory one optimizes over functions of only three variables, as opposed to functions of $3N$ variables in Schr\"odinger theory.

While Thomas--Fermi theory had been around and had been used for a long time, it was only in the 1970s that Lieb and Simon rigorously established its mathematical foundations \cite{tftheorysimon}. They answer, among other things, questions about the existence and uniqueness of solutions, their regularity and decay. Later, Lieb and collaborators investigated systematically density functional theories that are refinements of Thomas--Fermi theory. These findings are summarized in the review \cite{tftheoryrev}.

Another fundamental result that Lieb and Simon prove in their Thomas--Fermi paper is a rigorous relation between the ground state energy of the Hamiltonian $H$ in \eqref{eq:hamiltonian} and the minimal Thomas--Fermi energy in the joint limit $N\to \infty$ and $Z\to\infty$ with $N/Z\to\lambda\in (0,\infty)$. (Here we consider, for simplicity, the atomic case $K=1$ and set $Z=Z_1$. The most important case is $N=Z$, that is, the case of a neutral atom.) It is shown that
$$
\lim \frac{\inf\spec H}{\inf \mathcal E} = 1 \,.
$$
This convergence of energies is supplemented by convergence results of the one-particle densities of (approximate) ground states of the Schr\"odinger Hamiltonian. Technically, this result is related to semiclassical analysis, but outside of the typical regularity assumptions in this theory. The Lieb--Simon result has become a blueprint for other derivations of effective theories in scaling limits.

As an aside, we also mention the first proof, given by Lieb and Simon, of the existence of solutions to the Hartree--Fock equations for atoms and molecules \cite{hf}. Hartree--Fock theory is a more precise approximation to Schr\"odinger theory than Thomas--Fermi theory and used in the computation of atomic and molecular energies. In contrast to Thomas--Fermi theory, which is a density functional theory, where the unknown is a scalar function, Hartree--Fock is a density matrix theory, where the unknown is an operator. The Lieb--Simon paper is a foundational paper in the noncommutative calculus of variations.

In Thomas--Fermi theory, the Coulomb repulsion between the electrons is approximated by the quantity $D[\rho]$. In 1981, Lieb and Oxford \cite{lox} proved a lower bound on the difference between these two quantities, which became known as the Lieb--Oxford bound for the exchange energy. This bound is well known to quantum chemists and guides their thinking about the exchange correlation energy in molecules.

In 1983, Lieb published the paper ``Density Functionals for Coulomb Systems'' \cite{dft}, which laid out the theoretic foundations of density functional theory and is widely cited. It introduced a universal functional, known as the Levy--Lieb functional, which gives the lowest energy that can be reached with all possible quantum states having a given density function. This functional yields, by definition, the ground state energy of interacting quantum Coulomb systems (even if it is not known explicitly). This point of view has played a very important role. Density functional theory has exploded in the 1990s and is widely used in industry. It is now the method of choice to compute the quantum state of molecules and solids.

Among Lieb's more recent works in this direction are a mathematically rigorous justification of the Local Density Approximation in density functional theory, and a proof of the equivalence in the thermodynamic limit of three different definitions of the minimum energy of a homogeneous electron gas. These are joint works with Lewin and Seiringer \cite{unifelegas,lda}.


\subsection*{Lieb--Thirring inequalities}

Arguably the least obvious part in the approximation of the Schr\"odinger energy functional by the Thomas--Fermi functional is that the correct kinetic energy
$$
\int_{\R^3}\cdots \int_{\R^{3}} \sum_{n=1}^N |\nabla_n\psi|^2\,dx_1\ldots\,dx_N
$$
is replaced by the term
$$
\gamma^{\rm TF} \int_{\R^3} \rho(x)^\frac53\,dx
$$
for a certain explicit constant $\gamma^{\rm TF}>0$. It is in this step (and only in this step) of the approximation that the fermionic nature of the wave function $\psi$ enters. Behind this approximation is the observation that the kinetic energy per unit volume of a noninteracting Fermi gas in its ground state with constant density $\rho$ is $\gamma^{\rm TF} \rho^\frac53$. Imagining that particles described by a wave function $\psi$ with low energy are locally essentially in the ground state of the Fermi gas with the corresponding local density one arrives at the Thomas--Fermi expression for the kinetic energy.

The question arises whether this approximation can be substantiated by rigorous bounds. This is accomplished by the famous Lieb--Thirring inequality, which was a crucial ingredient in their proof of stability of matter \cite{som}. This inequality states that for any antisymmetric, normalized $\psi$ on $\R^{3N}$, one has
\begin{equation}
	\label{eq:lt}
	\int_{\R^3}\cdots \int_{\R^{3}} \sum_{n=1}^N |\nabla_n\psi|^2\,dx_1\ldots\,dx_N
	\geq K \int_{\R^3} \rho_\psi(x)^\frac53\,dx
\end{equation}
with
$$
\rho_\psi(x) = \sum_{n=1}^N \int_{\R^3}\cdots\int_{\R^3} |\psi(x_1,\ldots,x_{n-1},x,x_{n+1},\ldots,x_N)|^2\,dx_1\cdots dx_{n-1}\,dx_{n+1}\cdots dx_N \,.
$$
The important point here is that the constant $K$ is independent of $N$.

The Lieb--Thirring inequality \eqref{eq:lt} can be viewed as a mathematical expression of both the exclusion and the uncertainty principles in quantum mechanics. The connection with the exclusion principle is that the constant $K$ is independent of $N$ -- if symmetric functions $\psi$ (describing a bosonic system) would be allowed, the constant $K$ in \eqref{eq:lt} would have to deteriorate with $N$. (As an example, take $\psi$ as a product function $\phi(x_1)\cdots\phi(x_N)$.) For $N=1$ the Lieb--Thirring inequality reduces to a certain Sobolev interpolation inequality and for $N\geq 2$ it can be considered as a generalization thereof. The conclusion of Sobolev-type inequalities, namely that an $L^p$ with a ``large" $p$ can be controlled, is a non-concentration result and, therefore quantifies the uncertainty principle in quantum mechanics.

The constant $K$ in \eqref{eq:lt} that Lieb and Thirring obtained was smaller than $\gamma^{\rm TF}$, but it retained the important feature of being independent of $N$ (and of $\psi$, of course). The famous Lieb--Thirring conjecture states that the inequality should be valid with constant equal to $\gamma^{\rm TF}$. This would mean that the Thomas--Fermi approximation for the kinetic energy is a universal lower bound to its Schr\"odinger expression. There has been a lot of work on this constant, leading to the currently best bound of $(1.456)^{-2/3}\gamma^{\rm TF}$.

Lieb and Thirring did not prove inequality \eqref{eq:lt} directly, but first showed that it is equivalent to a certain inequality about sums of negative eigenvalues of one-body Schr\"odinger operators, and then verified the latter inequality. In their follow-up paper \cite{lt} they extended the latter inequality to arbitrary dimensions and arbitrary powers of eigenvalues. They proved that the negative eigenvalues $(E_j)$ of the Schr\"odinger operator $-\Delta+V$ in $L^2(\R^d)$ satisfy
$$
\sum_j |E_j|^\gamma \leq L_{\gamma,d} \int_{\R^d} V(x)_-^{\gamma+d/2}\,dx
$$
for all $\gamma>1/2$ if $d=1$ and $\gamma>0$ if $d\geq 2$. For $\gamma=1$ and $d=3$, this inequality is equivalent to \eqref{eq:lt}. Soon afterwards, Lieb \cite{clr} proved the corresponding inequality in the endpoint case, known as Cwikel--Lieb--Rozenblum inequality, namely
$$
\#\{ j:\ E_j <0 \} \leq L_{0,d} \int_{\R^d} V(x)_-^{d/2}\,dx
$$
for $d\geq 3$.

The general form of the famous Lieb--Thirring conjecture concerns the optimal values of the constants $L_{\gamma,d}$. Apart from the obvious relevance of the values of these constants in applications, the conjecture addresses on a conceptual level the strength of the exclusion principle in different dimensions. Lieb, together with Hundertmark and Thomas, have proved the only known case of an optimal inequality where the constant is not given by that arising from a semiclassical (or Thomas--Fermi-like) approximation. Also, after more than four decades, Lieb's value of the constant $L_{0,3}$ is still the smallest one that is known.

The Lieb--Thirring and Cwikel--Lieb--Rozenblum inequalities and their generalizations are of great importance in the study of large fermionic systems. Apart from that they have found applications in the context of nonlinear evolution equations like the Navier--Stokes equation \cite{navierstokes}. Moreover, they have attracted great interest from a purely mathematical point of view and the fact that orthogonality of functions leads to an improved dependence on the number of functions has been verified in a number of other functional inequalities as well, including \cite{bessel,strichartz}.


\subsection*{The ionization problem}

Let us return to the quantum many-body Hamiltonian $H$ in \eqref{eq:hamiltonian}. The infimum of its spectrum may or may not be an eigenvalue. If it is, we interpret the corresponding eigenfunction as describing the ground state of the system and we think of the $N$ electrons as bound to the nuclei. Physical intuition suggests that given nuclei with given charges can only bound a finite number of electrons, but even this is not quite obvious mathematically. The quantitative version of this question, namely how many electrons an atom (or a molecule) can bind, is still not settled, despite serious efforts.

For simplicity, let us restrict our attention to the atomic case, that is, $K=1$ in \eqref{eq:hamiltonian}. Experimental data and numerical estimates show that a nucleus of charge $Z$ can bind at most $Z+1$, or possibly $Z+2$ electrons. To prove (or disprove) this rigorously in the above Schr\"odinger model is a famous open problem.

One of the few nonasymptotic results in this direction is due to Lieb \cite{ion} and states that an atomic nucleus can not bind $2Z+1$ or more electrons. The factor $2$ in front of $Z$ looks ``too large" for large $Z$, but, for instance, for $Z=1$, the bound is optimal.  Only three decades later was Lieb's result improved for large $Z$.

A striking discovery by Benguria and Lieb \cite{benlie} is that the purported bound on the excess charge would not be true if the electrons were bosons, that is, if the antisymmetry requirement on admissible functions was replaced by the symmetry requirement
$\psi(\ldots,x_n,\ldots,x_m,\ldots)= \psi(\ldots,x_m,\ldots,x_n,\ldots)$ for $n\neq m$. They showed that there is a number $\lambda>1$ (numerically, $\lambda\approx 1.21$) such that a ``bosonic atom" can bind at least $(\lambda+o(1)) Z$ electrons as $Z\to\infty$. As a consequence of the Benguria--Lieb result, the Pauli principle (i.e.\ the antisymmetry requirement) needs to enter any possible proof. This excludes, in particular, any naive, purely electrostatic argument. 

Returning to the original, fermionic case, one can ask whether there is at least asymptotic neutrality in the sense that, as $Z\to\infty$, the number of electrons that can be bound is $Z+o(Z)$. That this is indeed the case was proved by Lieb together with Sigal, Simon and Thirring \cite{asympneut}. Other researchers obtained later quantitative bounds on the $o(Z)$ remainder, but showing that it is bounded seems to be out of reach of current techniques.


\subsection*{Bosonic Systems}

While our focus in this section was mostly on fermionic systems, Lieb has also made many fundamental contributions to the study of bosonic systems. Among those are the following:
\begin{enumerate}
	\item[(a)] The construction, with Liniger, of a model of a one-dimensional interacting Bose gas \cite{liebliniger}. This model has served as a prototype for later theoretical developments, and it is also realized experimentally.
	\item[(b)] Together with Yngvason, Lieb proved an asymptotic formula, conjectured 58 years earlier, for the ground state energy of a dilute Bose gas \cite{bec}. Subsequently, together with Seiringer and Yngvason, he rigorously derived the Gross--Pitaevskii equation for the ground state energy of dilute bosons in a trap, starting with many-body quantum mechanics \cite{gp}. This result has had a tremendous impact on the development of mathematical physics in the past two decades.
	\item[(c)] Together with Conlon and Yau \label{n75} and later with Solovej \cite{bogo1,bogo2}, Lieb proved the $N^{7/5}$ law for charged bosons. This means, roughly, that the energy does not obey a linear bound as in \eqref{eq:som}, but rather decreases like $-CN^{7/5}$. This was the first rigorous validation of Bogolubov's pairing theory of the Bose gas and paved the way for many current developments.
\end{enumerate}


\section{Functional Inequalities}

Lieb's name is inseparably connected with the subject of inequalities and a whole 700 page volume of his Selecta is dedicated to this topic \cite{selecta2}. In the previous section, in the discussion of the Lieb--Thirring inequality, we have already seen one instance of a functional inequality. In this section we will review three more examples, namely the entropy inequalities in matrix analysis, the Brascamp--Lieb inequalities and the sharp form of the Hardy--Littlewood--Sobolev inequality.

\subsection*{Lieb's concavity theorem and the strong subadditivity}

In 1973, Lieb and Ruskai proved the \emph{strong subadditivity of the quantum-mechanical entropy} \cite{ssa}. This theorem has many equivalent formulations, such as the concavity of the conditional entropy, the joint convexity of the relative entropy or the monotonicity of the relative entropy. Strong subadditivity, or one of its equivalents, has come to play an essential role in the modern and very active area of quantum information and quantum computation.

Let us state this theorem in its monotonicity formulation, originally derived by Lindblad in 1974 from a result of Lieb. A density matrix is a Hermitian, positive semi-definite matrix of trace one. The relative entropy (or Kullbach--Leibler divergence) of two density matrices $\rho$ and $\sigma$ is defined to be
$$
D(\rho\|\sigma) = \Tr \rho\ln\rho - \Tr \rho\ln\sigma \,.
$$
This quantity is nonnegative and vanishes if and only if $\rho=\sigma$. Roughly speaking, it measures how distinguishable $\rho$ and $\sigma$ are, even though the quantity is not symmetric in $\rho$ and $\sigma$. Quantum operations are described by completely positive, trace preserving maps, and the theorem of monotonicity of the relative entropy (also known as the \emph{data processing inequality}) states that for any such operation $\mathcal E$, one has
$$
D(\mathcal E\rho\|\mathcal E\sigma) \leq D(\rho\|\sigma) \,.
$$
Thus, applying a quantum operation can only make the states harder to distinguish. This makes it clear that the monotonicity of the relative entropy lies at the very foundation of quantum information theory.

Both the Lieb--Ruskai proof of the strong subadditivity of the entropy and the Lindblad proof of monotonicity of the relative entropy rest on a deep theorem that Lieb proved in his 1973 paper ``Convex trace functions and the Wigner–Yanase–Dyson conjecture" \cite{wy}. This theorem states that, for $p,q\geq 0$ with $p+q\leq 1$, the map $(A,B)\mapsto \Tr A^p B^q$, defined on nonnegative Hermitian matrices, is jointly concave.

Lieb's paper has generated considerable further work containing alternative proofs, generalizations and applications. He himself, often jointly with Carlen, has also returned several times to this circle of ideas and deepend our understanding of matrix analysis; see, e.g., \cite{hanner,minkowski}.


\subsection*{The Brascamp--Lieb inequalities}

In 1976 Brascamp and Lieb published the paper ``Best Constants in Young's Inequality, Its Converse and Its Generalization to More Than Three Functions" \cite{brascamp}. As the title suggests, this paper treats three related, but different topics and it is famous for all three.

The first one concerns a basic inequality in real analysis, namely the Young inequality
\begin{equation}
	\label{eq:y}
	\left| \iint_{\R^d\times\R^d} f(x)g(x-y) h(y) \,dx\,dy \right| \leq C_{p,q,r,d} \|f\|_{L^p(\R^d)} \|g\|_{L^q(\R^d)} \|h\|_{L^r(\R^d)}
\end{equation}
for three functions $f,g,h$ on $\R^d$ and parameters $1\leq p,q,r\leq\infty$ satisfying $\frac{1}{p} + \frac1q + \frac1r = 2$. This inequality appears frequently in theory and applications, since convolution is a basic operation in mathematical analysis and the Young inequality bounds its effect in Lebesgue spaces.

The question is to find the optimal (that is, smallest possible) constant $C_{p,q,r,d}$ in \eqref{eq:y} and to characterize all cases of equality. The first task was solved by Beckner around the same time as \cite{brascamp}. Brascamp and Lieb give an alternative proof and solve the second task. Their proof combines in an original way the technique of symmetric decreasing rearrangement with the operation of taking tensor products, leading to a convergence result in the spirit of the central limit theorem. It follows that the optimal constant is determined by Gaussian functions and attained only for those.

The second topic of the Brascamp--Lieb paper is that Young's inequality holds, with the inequality sign reversed, if $0<p,q\leq 1$ and $f,g,h$ are nonnegative. Again the authors are able to compute the optimal constant. This reverse Young's inequality contains, as a limiting case, the Pr\'ekopa--Leindler theorem.

The third topic of the Brascamp--Lieb paper is a far-reaching generalization of Young's inequality to an arbitrary number of functions, together with a replacement of the linear functions $x$, $y$ and $x-y$ on $\R^{2d}$, appearing in \eqref{eq:y}, by a much larger class of linear functions. The resulting family of inequalities is now known as Brascamp--Lieb inequalities. It contains not only Young's inequality, but also those of H\"older and of Loomis and Whitney as special cases. Brascamp and Lieb show that for given linear maps the inequality holds with a finite constant if and only if it holds for Gaussian functions, and in this case the optimal constant can be computed using the latter class of functions. Their argument is again based on a central limit theorem.

Besides the original application to statistical mechanics, the Brascamp--Lieb inequalities have come to play an important role in convexity theory, harmonic analysis and other parts of mathematics.

In his 1990 paper with the title ``Gaussian Kernels Have Only Gaussian Maximizers" \cite{gaussian} Lieb revisits the topic of Young and Brascamp--Lieb inequalities. He proves a general theorem about the norm of a large class of integral operators. Among other things, he characterizes the cases of equality in the Hausdorff--Young inequality with optimal constant. Again, they are given by Gaussian functions.

We finally mention that there is yet another celebrated inequality known as Bras\-camp--Lieb inequality. This is a Poincar\'e (or spectral gap) inequality for log-concave probability distributions \cite{brascamp2}.

\subsection*{The sharp Hardy--Littlewood--Sobolev inequality}

A fundamental result of real analysis, known as the Hardy--Littlewood--So\-bo\-lev inequality, the weak Young inequality or the theorem of fractional integration states that, for $0<\lambda<d$ and $1<p,r<\infty$ with $\frac1p + \frac{\lambda}{d} + \frac 1r = 2$, one has
\begin{equation}
	\label{eq:hls}
	\left| \iint_{\R^d\times\R^d} \frac{f(x)\,h(y)}{|x-y|^\lambda}\,dx\,dy \right| \leq C_{\lambda,p,r} \|f\|_{L^p(\R^d)} \|h\|_{L^r(\R^d)} \,.
\end{equation}
This inequality has many applications in pure and applied mathematics. For instance, in the case $\lambda=d-2$, $d\geq 3$, the kernel in this inequality is the Coulomb kernel and the quantity on the left side represents the Coulomb interaction of two charge densities $f$ and $g$. For $\lambda=d-2$ and $p=r$, inequality \eqref{eq:hls} is equivalent to the Sobolev inequality.

Lieb's 1983 paper ``Sharp constants in the Hardy--Littlewood--Sobolev and related inequalities'' \cite{hls} had a profound impact on the field of Calculus of Variations. Indeed, this paper and his related 1983 paper ``A relation between pointwise convergence of functions and convergence of functionals'' \cite{brezislieb} with Brezis are Lieb's most cited pure mathematics papers.

The ``Sharp constants" paper is remarkable for at least two distinct aspects, namely for the development of a rather general compactness argument and for an ingenious observation concerning the problem at hand. We briefly discuss these two points.

The compactness aspect of Lieb's HLS paper concerns the question whether there is a nontrivial pair $(f,h)$ such that equality holds in \eqref{eq:hls} with $C_{\lambda,p,r}$ taken to be the minimal constant. Lieb had worked earlier on questions about existence of optimizers for certain variational problems. In contrast to classical variational problems, these problems were often translation-invariant and so Lieb had to find methods to deal with the corresponding loss of compactness (for instance, by symmetric decreasing rearrangement \cite{choquard} or by an original ``running-after argument'' \cite{intersection}). The optimization problem for the HLS inequality, however, features another potential loss of compactness, namely through dilations. In order to deal with these problems, Lieb found a strengthening of Fatou's lemma, which says that, if functions $f_j$ on a measure space $X$ are bounded in $L^p$ and converge pointwise a.e.\ to a function $f$, then
\begin{equation}
	\label{eq:bl}
	\lim_{j\to\infty} \int_X \left| |f_j|^p - |f|^p - |f_j-f|^p \right|dx = 0 \,.
\end{equation}
In particular,
$$
\int_X |f_j|^p\,dx = \int_X |f|^p\,dx + \int_X |f_j-f|^p\,dx + o(1) \,.
$$
For comparison, in Fatou's lemma, the second term on the right side is omitted, thus leading to an inequality rather than an equality. The statement \eqref{eq:bl}, with $|\cdot|^p$ replaced by more general functions, is known as the Brezis--Lieb lemma and constitutes a fundamental tool in functional analysis and the calculus of variations.

Lieb's compactness method, sometimes called the method of the missing mass, tracks carefully the remainder term $f_j-f$ in \eqref{eq:bl}. This technique is quite robust and has found various applications to different settings, including \cite{steintomas}. Also other compactness methods use the Brezis--Lieb lemma as a fundamental ingredient.

The second remarkable aspect of Lieb's HLS paper concerns the special case $p=r$ (but $0<\lambda<d$ is arbitrary). Lieb managed to compute the smallest possible constant $C_{\lambda,p,r}$ explicitly and to characterize all pairs $(f,h)$ of functions for which equality is attained. The crucial observation in Lieb's proof was a ``hidden" symmetry, namely the conformal invariance of \eqref{eq:hls} for $p=r$. Lieb's paper has spawn a field of variational problems with conformal invariance. Among these developments is also the sharp form of the HLS inequality on the Heisenberg group \cite{heisenberg}.


\section{Topics not covered}

In the previous sections we have discussed two topics of Lieb's work, namely quantum Coulomb systems and functional inequalities. Those make up just a fraction of the body of Lieb's work and it is not unlikely that other writers would have picked completely different topics. We would be remiss if we would not at least briefly mention a few more.

We have completely ignored the important chapter of exactly soluble models in Lieb's work in the 1960s. This field was revolutionized by Lieb's solutions of the ice problem, the Fierz--Rys F-model and Slater's KDP model, as well as the already mentioned Lieb--Liniger model. In connection with ice type models, Lieb and Temperley constructed what became known as the Temperley--Lieb algebra, which has applications in several areas of mathematics, e.g., knot theory, and in physics. Lieb's work lies at the foundations of modern statistical mechanics and has a long-lasting influence on integrable probability and combinatorics, to name just a few.

Concerning Lieb's contributions to condensed matter physics Nachtergaele, Solovej and Yngvason write in the preface of the third volume of Lieb's selected works \cite{selecta3}:
\begin{quote}
	``The impact of Lieb's work in mathematical condensed matter physics is unrivaled. It is fair to say that if one were to name a founding father of the field, Elliott Lieb would be the only candidate to claim this singular position."
\end{quote}
This area includes, in particular, Lieb's seminal work on the Hubbard model, including its solution in the one-dimensional case with Wu and the discovery of Lieb ferromagnetism and the Lieb lattice, as well as the highly cited joint works with Schultz and Mattis on two soluble models of an antiferromagnetic chain and with Mattis on the Luttinger model and the discovery of bosonization.

Already in the introduction we mentioned the ``Lieb--Robinson bounds”, which establish the deep fact that there is a finite group velocity for information propagation in quantum spin systems. This turned out to be quite important for ideas about quantum computation and condensed matter physics.

Moreover, together with Affleck, Kennedy and Tasaki, Lieb developed what came to be known as the AKLT model of a spin-one spin chain. This has proved to be an important prototype of a class of models and eigenstates known as matrix product states.

From the area of statistical mechanics we mention Lieb's proof of the existence of phase transitions in classical and quantum spin systems, obtained jointly with Fr\"ohlich, Israel and Simon and with Dyson and Simon. The proofs are based on the method of reflection positivity, which Lieb has masterfully used in several situations.

Another famous result is Lieb's work with Heilman on the zeros of the partition function of the monomer-dimer problem, which is related to the matching problem in combinatorics. This is also of importance in computer science.

A notable mention is a new approach to the physics and mathematics of the second law of thermodynamics, developed jointly with Yngvason.

The list of topics in this section has so far focused on the more physics-oriented work of Lieb. From the more mathematics-oriented work let us just mention his fundamental contributions to the theory of symmetric decreasing rearrangement. This includes, among other things, a general rearrangement inequality for many functions, obtained with Brascamp and Luttinger, as well as the proof that symmetric decreasing rearrangement can be discontinuous in Sobolev spaces, obtained with Almgren.

Finally, we would like mention the book by Lieb and Loss \cite{LiLo}, which has become a standard textbook for graduate courses in mathematical analysis and which eloquently and concisely promotes the philosophy of rigorous mathematics with a view towards applications.

\medskip

We hope that these pages may serve as an invitation to look at Lieb's original research papers. Only in this way can the reader feel the clarity, vitality and beauty of Lieb's work, a work that has, and continues to, inspire generations.

\medskip

Congratulations on receiving the Gauss Prize, Elliott!


\bibliographystyle{amsalpha}

\end{document}